\documentclass{amsart}
\usepackage{epsfig}
\usepackage{amsfonts,amssymb,amscd,amsmath,euscript,enumerate,verbatim,calc}
\theoremstyle{plain}

\newtheorem{theorem}{Theorem}[section]
\newtheorem{lemma}[theorem]{Lemma}
\newtheorem{proposition}[theorem]{Proposition}

\newtheorem{conjecture}[theorem]{Conjecture}

\theoremstyle{definition}
\newtheorem{definition}[theorem]{Definition}

\newtheorem{remark}[theorem]{Remark}
\newtheorem{example}[theorem]{Example}
\def\F{{\mathbb F}}
\def\Fq{{\mathbb F}_q}
\def\Fqm{{\mathbb F}_{q^m}}
\def\Fqmn{{\mathbb F}_{q^{mn}}}

\newcommand{\ord}{\operatorname{ord}}
\newcommand{\GL}{\operatorname{GL}}

\newcommand{\im}{\operatorname{im}}
\newcommand{\Cf}{{\sf C}_f}
\newcommand{\bcm}{{\rm BCM}(m,n;q)}
\newcommand{\bcms}{{\rm BCMS}(m,n;q)}
\newcommand{\bcmsone}{{\rm BCMS}(m,1;q)}
\newcommand{\pmn}{{\EuScript P}(mn;q)}
\newcommand{\pmone}{{\EuScript P}(m;q)}
\begin{document}
\title[Primitive Polynomials, Singer Cycles and $\sigma$-LFSR\lowercase{s}]{Primitive Polynomials, Singer Cycles, and Word-Oriented Linear Feedback Shift Registers}
\author{Sudhir R. Ghorpade}
\address{Department of Mathematics, 
Indian Institute of Technology Bombay,\newline \indent
Powai, Mumbai 400076, India.}
\email{srg@math.iitb.ac.in}

\author{Sartaj Ul Hasan}
\address{Department of Mathematics,
Indian Institute of Technology Bombay,\newline \indent
Powai, Mumbai 400076, India \newline \indent
and \newline \indent
Scientific Analysis Group, Defense Research and Development Organisation  \newline \indent
Metcalfe House, Delhi 110054, India}
\email{sartajulhasan@gmail.com}

\author{Meena Kumari}
\address{Scientific Analysis Group, Defense Research and Development Organisation  \newline \indent
Metcalfe House, Delhi 110054, India}
\email{rameena10@yahoo.co.in}


\begin{abstract}
Using the structure of Singer cycles in general linear groups, we prove that a 
conjecture of Zeng, Han and He (2007) holds in the affirmative in a special case, and
outline a plausible approach to prove it in the general case. 
This conjecture is about the number of primitive $\sigma$-LFSRs of a given order over a finite field, and it generalizes a known formula for the number of primitive LFSRs, which, in turn, is the 
 number of primitive polynomials of a given degree over a finite field. Moreover, this conjecture is intimately related to an open question of Niederreiter (1995) on the enumeration 
of splitting subspaces of a given dimension. 
\end{abstract}
\date{\today}
\maketitle
\section{Introduction}

Denote, as usual, by $\Fq$ the finite field with $q$ elements and by $\Fq[X]$ the ring of polynomials in one variable $X$ with coefficients in $\Fq$.
It is elementary and well known that if $f(X)\in \Fq[X]$ is of degree $n$ and $f(0)\ne 0$, then 
$f(X)$ divides $X^e-1$ for some positive integer $e\le q^n-1$. The least such $e$ is called the \emph{order} of $f(X)$ 
and is denoted by $\ord f(X)$. We say that a monic polynomial $f(X)\in \Fq[X]$ of degree $n$ is \emph{primitive} if $f(0)\ne 0$ and $\ord f(X) = q^n-1$. The study of primitive polynomials goes back to Gauss and is an interesting and important part of the theory of finite fields. A basic reference is \cite[Ch. 3]{LN} and 
some of the relevant facts about primitive polynomials 
are stated in Section \ref{sec1} below. 

Elements of the maximum possible order in the finite group $\GL_n(\Fq)$ of $n\times n$ nonsingular matrices with entries in $\Fq$ are called \emph{Singer cycles}. These are closely related to primitive polynomials since this maximum possible order is, in fact, $q^n-1$, and moreover, characteristic polynomials of Singer cycles are primitive. 
We refer to \cite{HUP} and \cite{SING} for some basic aspects of the study of Singer cycles and provide, for the convenience of the reader, a brief outline of basic results together with consequences that are useful for this paper in Section \ref{sec2}.

Linear feedback shift registers (LFSRs) are devices frequently used in cryptography and coding theory (cf. [\ref{GOL}, \ref{bib13}]). In effect, a LFSR can be viewed as an infinite sequence of elements of $\Fq$ generated by
finitely many initial values and a homogeneous linear recurrence relation. In the binary case ($q=2$), these sequences are used for efficient encryption of data in designing stream ciphers. In general, it can be shown that these sequences are (ultimately) periodic and the maximum possible period of an $n^{\rm th}$ order linear recurring sequence is $q^n-1$. 
(See Section \ref{sec1} for details.)  
In order to have good cryptographic properties [\ref{GOL}], one is mainly interested in the sequences that  have the maximum period. The LFSRs corresponding to sequences with maximum period are known as \emph{primitive LFSRs}. 
Using the connection with primitive polynomials or otherwise, it 
is readily seen that the number of primitive LFSRs of order $n$ over $\Fq$ 
is given by 
\begin{equation}
\label{NoLFSR}
	\frac {\phi(q^{n}-1)}{n}
\end{equation}
where $\phi$ is the Euler totient function. 

In this paper, we consider a recent generalization due to Zeng, Han and He [\ref{bib23}] of a (traditional) LFSR to a word-oriented linear feedback shift register, called $\sigma$-LFSR. It is argued in [\ref{bib23}] that the $\sigma$-LFSRs meet the dual demands of high efficiency and good cryptographic properties, and that these can be viewed as a solution to a problem of Preneel \cite{BP} on designing fast and secure LFSRs with the help of the
word operations of modern processors and the techniques of parallelism. Notions of primitivity readily extend from LFSRs to $\sigma$-LFSRs although the connection with primitive polynomials and matrices is a little more intricate. 
Unlike \eqref{NoLFSR}, a simple formula for the number of primitive $\sigma$-LFSRs of order $n$ over $\Fqm$ is not known, but an intriguing explicit formula in the binary case has been conjectured. 
The main aim of this paper is to elucidate, extend and understand this conjectural formula of 
Zeng, Han and He [\ref{bib23}]. In general, the conjecture is that the number of primitive $\sigma$-LFSRs of order $n$ over $\Fqm$ is given by
\begin{equation}
\label{NoSigmaLFSR}
\frac{\phi(q^{mn}-1)}{mn}.q^{m(m-1)(n-1)} \displaystyle \prod_{i=1}^{m-1}(q^m-q^i).
\end{equation}

After a preliminary version of this paper was prepared, we found that the seemingly new notion of a $\sigma$-LFSR can, in fact, be traced back 
to the work of Niederreiter (1993-1996) mainly in the context of pseudorandom number generation. Indeed, in a series of papers [\ref{Niederreiter1}, \ref{Niederreiter2}, \ref{Niederreiter3}, \ref{Niederreiter4}], Niederreiter has introduced the so called \emph{multiple recursive matrix method} and the notion of recursive vector sequences. The latter are essentially the same as sequences generated by a $\sigma$-LFSR, modulo a natural isomorphism between 
the field $\Fqm$ with $q^m$ elements and the vector space $\Fq^m$ 
of dimension $m$ over 
The question of counting the number of primitive $\sigma$-LFSRs of a given order $n$ over $\Fqm$ is considered in [\ref{Niederreiter2}, p. 11] under a different guise (cf. Remark~\ref{rem:Nied}), and is termed as open problem. 
However, no explicit formula for this number is given, even conjecturally, in the work of Niederreiter, and therefore, the credit for formulating \eqref{NoSigmaLFSR} should go to Zeng, Han and He [\ref{bib23}] at least in the binary case. Moreover, in a personal communication, 
Professor Niederreiter has informed us that the problem of counting the number of primitive $\sigma$-LFSRs of a given order $n$ over $\Fqm$ 
is  still open to the best of his knowledge. 
%

Our main results are as follows. 
We work throughout in the general $q$-ary case and first give an alternative formulation of the conjecture in terms of the enumeration of certain Singer cycles (Theorem \ref{conj2}). Next, we give a plausible approach to derive \eqref{NoSigmaLFSR} by noting that it suffices to analyze the image and the fibers of a natural map from a certain class of $mn\times mn$ matrices to the set of primitive polynomials of degree $mn$. We accomplish the first task by showing that this map is surjective (Theorem \ref{surjprop}). As for the second, we give a conjectural description of the fibers  (Conjecture \ref{twoconj}).  Moreoever, we use certain properties of Singer cycles to prove that \eqref{NoSigmaLFSR} as well as the more refined Fiber Conjecture hold when $n=1$ and $m$ is arbitrary (Theorem~\ref{onecase}). It may be noted that in the 
other initial case $m=1$, \eqref{NoSigmaLFSR} is an immediate consequence of \eqref{NoLFSR}. 

This paper is written in a fairly self-contained manner with the hope that it would stimulate some interest even among those that are not interested in cryptographic applications {\it per se}, in proving formula \eqref{NoSigmaLFSR} and taking up allied problems.

\section{Primitive Polynomials and Primitive LFSRs}
\label{sec1}

By a \emph{primitive} element in a finite cyclic group $G$ we mean a generator of $G$. 
Primitive polynomials in $\Fq[X]$, as defined in the Introduction, are related to primitive elements 
by the following characterization \cite[Thm. 3.16]{LN}, which is sometimes used to give an alternative definition of primitive polynomials. 

\begin{proposition}
\label{PrimPolyElt}
Let $f(X)\in \Fq[X]$ be of degree $n \ge 1$. Then $f(X)$ is a primitive polynomial if and only if $f(X)$ is the minimal polynomial of a primitive element of the cyclic group $\mathbb F_{q^n}^*$ of nonzero elements of the finite field $\mathbb F_{q^n}$. 
\end{proposition}

Using the above theorem together with the fact that the 
number of primitive elements in a cyclic group of order $N$ is $\phi (N)$, 
we readily see that the number of primitive polynomials in $\Fq[X]$ of degree $n$ is given by \eqref{NoLFSR}. 

We shall now proceed to review the basic definitions and some of the basic results concerning 
linear feedback shift registers.
     
\begin{definition}
\label{def:lfsr}
Let $n$ be a positive integer and let $c_0,c_1, \dots, c_{n-1} \in \mathbb F_q$. Given any $n$-tuple $(s_0,s_1, \dots, s_{n-1}) \in \mathbb F_{q}^n$, let $s^{\infty}=(s_0,s_1,\dots)$ denote the infinite sequence of elements of $\Fq$  determined by the following linear recurrence relation:
\begin{eqnarray}
s_{i+n}=s_ic_0+s_{i+1}c_1+\cdots +s_{i+n-1}c_{n-1} \quad \mbox{for} \quad i=0,1,\dots \label{lfsr} 
\end{eqnarray}
The system (\ref {lfsr}) is called a \emph{linear feedback shift register (LFSR)} of order $n$ over $\Fq$, while the sequence $s^{\infty}$ is referred to as the \emph{sequence generated by the LFSR} (\ref{lfsr}). The $n$-tuple $(s_0,s_1, \cdots, s_{n-1})$ is called the \emph{initial state} of the LFSR (\ref{lfsr}) and the polynomial $X^n -c_{n-1}X^{n-1}- \cdots  -c_1X-c_0 $ is called the \emph{characteristic polynomial} of the LFSR (\ref{lfsr}). 
The sequence $s^{\infty}$ is said to be 
\emph{ ultimately periodic} if there are integers $r, n_0$ with $r\ge 1$ and $n_0\geq0$ such that 
$s_{j+r}=s_j$ for all $j \geq n_0$. The least positive integer $r$ with this property is called the \emph{period} of $s^{\infty}$ and the corresponding least nonnegative integer $n_0$ is called the \emph{preperiod} of $s^{\infty}$. 
The sequence  $s^{\infty}$ is said to be \emph{periodic} if its preperiod is $0$.
\end{definition}

Some basic facts about LFSRs are summarized in the two propositions below.  Proofs can be found, for example,  in \cite[Ch. 8]{LN}. 

\begin{proposition}
\label{ulti}  
For the sequence $s^{\infty}$ generated by the LFSR (\ref{lfsr}) of order $n$ over $\mathbb F_{q}$,
we have the following.
\begin{enumerate}
	\item[{\rm (i)}] $s^{\infty}$ is ultimately periodic and its period is $\le q^n-1$. 
	\item[{\rm (ii)}] 
	If $c_0\ne 0$, 
	then $s^{\infty}$ is periodic. Conversely, if $s^{\infty}$ is periodic whenever the initial state is of the form $(b,0, \dots, 0)$, where $b\in \Fq$ with $b\ne 0$, then $c_0\ne 0$. 
\end{enumerate} 
\end{proposition}

We say that a LFSR  of order $n$ over $\mathbb F_{q}$ is \emph{primitive} if for any choice of a nonzero initial state, the sequence generated by that LFSR is periodic of period $q^{n}-1$. 
Primitive LFSRs admit the following characterization. 

\begin{proposition}
\label{primlfsr}
A LFSR of order $n$ over $\mathbb F_{q}$ is primitive if and only if its characteristic polynomial is a primitive polynomial of degree $n$ in $\Fq[X]$.
\end{proposition}

As an immediate consequence of Propositions \ref{PrimPolyElt} and \ref{primlfsr}, we see that the number of primitive LFSRs of order $n$ over $\mathbb F_{q}$ is given by 
\eqref{NoLFSR}. 

\section{Singer Cycles and Singer Subgroups} 
\label{sec2}

The following result about orders of elements in a general linear group over finite field is well known. We include a more elaborate version and a quick proof since it seems a bit difficult to locate in or extract from the literature. 
An alternative (and somewhat longer) proof of the inequality below can be found, for example, in \cite[p. 742]{Dara}. 
In what follows, for an element $A$ of a finite group $G$, we denote by $o(A)$ the order of $A$ in $G$.
%
%

\begin{proposition}
\label{maxorder}
Let $A\in \GL_n(\Fq)$ and let $p(X)\in \Fq[X]$ be the minimal polynomial of $A$ and $\chi(X)\in \Fq[X]$ be the characteristic polynomial of $A$. 
Then $p(0)\ne 0$ and $o(A) = \ord p(X)$. In particular, 
$
o(A)\le q^n-1,
$  
and moreover, if the equality holds, then $p(X)= \chi (X)$. Also, we have: 
\begin{equation}
\label{equivalence}
o(A) = q^n-1 \;  \Longleftrightarrow  \; p(X) \mbox{ is primitive of degree $n$ } \Longleftrightarrow \; \chi (X) \mbox{ is primitive.}
\end{equation}
\end{proposition}

\begin{proof}
Since $A$ is nonsingular, $0$ is not an eigenvalue of $A$ and hence $p(0)\ne 0$. Now, if $I$ denotes the $n\times n$ identity matrix over $\Fq$, then for any positive integer $e$, we clearly have 
$$
A^e = I \; \Longleftrightarrow \; p(X) \mbox{ divides } X^e - 1. 
$$
Consequently, $o(A) = \ord p(X)$. Further, $\deg \chi (X)=n$ and in view of the Cayley-Hamilton Theorem, $p(X)$ divides 
$\chi(X)$. In particular, $\deg p(X) \le n$ and hence $o(A) = \ord p(X)\le q^n-1$.  Moreover, if $\ord p(X) = q^n-1$, then $\deg p(X)=n=\deg \chi (X)$, and hence $p(X)=\chi(X)$. On the other hand, if $\chi (X)$ is primitive, then it is irreducible and so $\chi(X)=p(X)$. 
This yields the equivalence in \eqref{equivalence}. 
\end{proof}

A cyclic subgroup of $\GL_n(\Fq)$ of order $e=  q^n-1$ is called a \emph{Singer subgroup} of $\GL_n(\Fq)$ and an element of $\GL_n(\Fq)$ of order $e$ is called a \emph{Singer cycle} in $\GL_n(\Fq)$. This terminology stems from \cite{SING} and seems appropriate since $\GL_n(\Fq)$ can be viewed as a subgroup of the symmetric group ${\mathfrak{S}}_e$ via the natural transitive action of $\GL_n(\Fq)$ on the set $\Fq^n\setminus \{0\}$, and elements of $\GL_n(\Fq)$ of order $e$ evidently correspond to $e$-cycles in ${\mathfrak{S}}_e$. We now recall two results from  \cite[II.\S 7]{HUP} (see also [\ref{COSSI}]) about Singer subgroups that will be useful to us later. 

\begin{proposition}
\label{conjugateSinger}
Any two Singer subgroups in  $GL_n(\mathbb F_q)$ are conjugate. 
\end{proposition}

\begin{proposition}
\label{normalizer}
Let $\sigma$ be the Frobenius automorphism of order $n$ of the field $\mathbb F_{q^n}$. Identify $\mathbb F_{q^n}$ 
with the vector space $\Fq^n$ and regard $\sigma$ as an element of $\GL_n(\Fq)$. Also, let $H$ be a Singer subgroup of $\GL_n(\Fq)$ and $N$ denote its normalizer in $\GL_n(\Fq)$. Then $N$ is isomorphic to the  semi-direct product $H \rtimes \left\langle \sigma\right\rangle$ of $H$ and the cyclic subgroup  of $\GL_n(\Fq)$ generated by $\sigma$. 
\end{proposition}

It may be noted that Proposition \ref{maxorder} relates Singer cycles to primitive polynomials. To work in the other direction, we can use companion matrices. Recall that if 
$f(X)= X^n -c_{n-1}X^{n-1}- \cdots  -c_1X-c_0 $ is a monic polynomial of degree $n\ge 1$ in $\Fq[X]$, then the
\emph{companion matrix} $\Cf$ of $f(X)$ is 
the $n \times n$ matrix 
$$ 
\Cf =
\begin {pmatrix}
0 & 0 & 0 & . & . & 0 & 0 & c_0\\
1 & 0 & 0 & . & . & 0 & 0 & c_1\\
. & . & . & . & . & . & . & .\\
. & . & . & . & . & . & . & .\\
0 & 0 & 0 & . & . & 1 & 0 & c_{n-2}\\
0 & 0 & 0 & . & . & 0 & 1 & c_{n-1}
\end {pmatrix} . 
$$
It is clear that $\det\Cf = (-1)^{n+1}c_0$. In particular, $\Cf\in \GL_n(\Fq)$ if and only if $f(0)\ne 0$. Also, we know from linear algebra that 
$f(X)$ is the minimal polynomial as well as the characteristic polynomial of $\Cf$. 
Thus, in view of Proposition \ref{maxorder}, we see that if $f(0)\ne 0$, then $\ord f(X) =  o(\Cf)$ and  
that $f(X)$ is a primitive polynomial if and only if $\Cf$ is a Singer cycle in $\GL_n(\Fq)$. In turn, primitive LFSRs of order $n$ over $\Fq$ are related to Singer cycles in $\GL_n(\Fq)$. To see the latter in a more direct way, it may be useful to observe that the companion matrix, say $A$, 
of the characteristic polynomial of the LFSR (\ref{lfsr}) is its
state transition matrix. 
Indeed, 
the $k^{\rm th}$ state $S_k:=\left(s_{k}, s_{k+1}, \dots , s_{k+n-1}\right)$ of the LFSR 
(\ref{lfsr}) is obtained from the initial state $S_0:=\left(s_{0}, s_{1}, \dots , s_{n-1}\right)$ by 
$S_k = S_0 A^k$,  for any $k\ge 0$.

\section {Word-Oriented Feedback Shift Register: $\sigma$-LFSR}
\label{seclfsr}

Given 
any ring $R$ and any positive integer $d$, let $M_d(R)$ denote the set of all $d\times d$ matrices with entries in $R$.
Fix throughout this and the subsequent sections, positive integers $m$ and $n$, and a vector space basis 
$\{\alpha_0, \dots, \alpha_{m-1}\}$ of ${\mathbb F}_{q^m}$  over $\mathbb F_q$. Given any $s\in {\mathbb F}_{q^m}$, there are unique $a_0, \dots, a_{m-1} \in {\mathbb F}_{q}$ such that 
$s = a_0 \alpha_0 + \cdots + a_{m-1}\alpha_{m-1}$, 
and we shall denote the corresponding co-ordinate vector 
$(a_0, \dots, a_{m-1})$ of $s$ by $\mathbf{s}$.  Evidently, the association $s\longmapsto \mathbf{s}$ gives a vector space isomorphism of $\mathbb F_{q^m}$ onto $\mathbb F_q^m$. Elements of $\mathbb F_q^m$ may be thought of as row vectors and 
so $\, \mathbf{s}C$ is a well-defined element of $\mathbb F_q^m$ for any 
$\mathbf{s} \in \mathbb F_q^m$ and $C\in M_m(\Fq)$. 
Following \cite{Zeng}, and in analogy with LFSRs, we define a ($q$-ary) $\sigma$-LFSR as follows. 

\begin{definition}
\label{def:sigmalfsr}
Let $C_0, C_1, \dots, C_{n-1} \in M_m(\mathbb F_q)$. Given any $n$-tuple $(s_0, \dots, s_{n-1})$  of elements of  $\mathbb F_{q^m}$, let $s^{\infty}=(s_0,s_1,\dots)$ denote the infinite sequence of elements of ${\mathbb F}_{q^m}$  
determined by the following linear recurrence relation:
\begin{eqnarray}
{\mathbf{s}}_{i+n}={\mathbf{s}}_iC_0+{\mathbf{s}}_{i+1}C_1+\cdots +{\mathbf{s}}_{i+n-1}C_{n-1} \quad \mbox{for} \quad i=0,1,\dots \label{sigmalfsr} 
\end{eqnarray}
The system (\ref {sigmalfsr}) is called a \emph{sigma linear feedback shift register ($\sigma$-LFSR)} of order $n$ over $\mathbb F_{q^m}$, while the sequence $s^{\infty}$ is referred to as the \emph{sequence generated by the $\sigma$-LFSR} (\ref{sigmalfsr}). The $n$-tuple $(s_0,s_1, \cdots, s_{n-1})$ is called \emph{initial state} of the $\sigma$-LFSR (\ref{sigmalfsr}) and the polynomial $X^n -C_{n-1}X^{n-1}- \cdots  -C_1X-C_0$ with matrix coefficients is called the \emph{$\sigma$-polynomial} of the $\sigma$-LFSR (\ref{sigmalfsr}). 
The sequence $s^{\infty}$ is said to be 
\emph{ ultimately periodic} if there are integers $r, n_0$ with $r\ge 1$ and $n_0\geq0$ such that 
$s_{j+r}=s_j$ for all $j \geq n_0$. The least positive integer $r$ with this property is called the \emph{period} of $s^{\infty}$ and the corresponding least nonnegative integer $n_0$ is called the \emph{preperiod} of $s^{\infty}$. 
The sequence  $s^{\infty}$ is said to be \emph{periodic} if its preperiod is $0$.
\end{definition}

%
%
%

The following analogue of 
Proposition \ref{ulti} is easily proved in a similar manner as in the classical case of LFSRs. 

\begin{proposition}
\label{sigmaulti}
For the sequence $s^{\infty}$ generated by the $\sigma$-LFSR (\ref{sigmalfsr}) of order $n$ over $\mathbb F_{q^m}$, 
we have the following. 
\begin{enumerate}
	\item[{\rm (i)}] $s^{\infty}$ is ultimately periodic, and its period is $\leq q^{mn}-1$. 
	\item[{\rm (ii)}] 
	If $C_0$ is nonsingular, then $s^{\infty}$ is periodic.
	Conversely, if $s^{\infty}$ is periodic whenever the initial state is of the form $(b, 0, \dots , 0)$, where $b\in \mathbb F_{q^m}$ with $b\ne 0$, then $C_0$ is nonsingular. 
\end{enumerate} 
\end{proposition}
 
We say that a $\sigma$-LFSR  of order $n$ over $\mathbb F_{q^m}$ is \emph{primitive} if for any choice of nonzero initial state, the sequence generated by that $\sigma$-LFSR is periodic of period  $q^{mn}-1$. 
In view of Proposition~\ref{sigmaulti}, if $X^n -C_{n-1}X^{n-1}- \cdots  -C_1X-C_0 \in M_m(\Fq)[X]$ is the $\sigma$-polynomial of a primitive $\sigma$-LFSR, then the matrix $C_0$  is necessarily nonsingular. 

Since the $\sigma$-polynomial of a $\sigma$-LFSR has coefficients in the noncommutative ring of matrices, notions such as irreducibility or primitivity are not readily applicable to it, and an analogue of Proposition \ref{primlfsr} is not obvious. 
However, as stated in \cite[Thm. 2]{Zeng} and proved in \cite[Thm. 3]{Zeng4} (see also \cite[Thm. 4]{N1}), we have the following 
characterization of primitive $\sigma$-LFSRs. 

\begin{proposition}
\label{primsigmalfsr}
Let $f(X)=X^n - C_{n-1}X^{n-1} - \cdots - C_1X - C_0 \in M_m(\mathbb F_q)[X]$ be the $\sigma$-polynomial of a 
$\sigma$-LFSR of order $n$ over $\mathbb F_{q^m}$, where $C_0 \in GL_m(\mathbb F_q)$ and $C_{\ell}\in M_m(\Fq)$ for $\ell =1, \dots , n-1$. 
For $1\le i, j \le m$, let 
$ f^{ij}(X)\in \Fq[X]$ be the polynomial of degree $n$ given by
$$
f^{ij}(X)=\delta^{ij}X^n - \sum_{\ell=0}^{n-1}c_{\ell}^{ij}X^{\ell}, 
$$             
where $\delta^{ij}$ is the Kronecker delta and $c_{\ell}^{ij}$ is  the $(i,j)^{\rm th}$ entry of the $m\times m$ matrix $C_{\ell}$ for \mbox{$\ell=0,1, \dots, n-1$.}  
Finally, let $\Delta(X)$ denote the determinant of the $m\times m$ matrix $\left(f^{ij}(X) \right)$ with polynomial entries. Then the $\sigma$-LFSR is primitive if and only if the $\Delta (X)$ is a primitive polynomial over 
$\mathbb F_q$ of degree $mn$.
\end{proposition}

The $q$-ary version of Conjecture 1 of \cite{Zeng} is the following. 

\begin{conjecture}
\label{conj1}
The number of primitive $\sigma$-LFSR of order $n$ over $\mathbb F_{q^m}$ is given by the formula \eqref{NoSigmaLFSR} stated in the Introduction. 
\end{conjecture}

We note that since $|GL_m(\mathbb F_q)| = (q^m-1) (q^m-q)\cdots (q^m - q^{m-1})$, the formula \eqref{NoSigmaLFSR} can be equivalently written as 
\begin{equation}
\label{NoSigmaLFSR2}
\Upsilon(m,n;q)= 
\frac{|GL_m(\mathbb F_q)|}{q^m-1}. \frac{\phi(q^{mn}-1)}{mn}.q^{m(m-1)(n-1)} 
\end{equation}
In fact, it appears in \cite{Zeng} in this form in the case $q=2$. As noted in \cite{Zeng}, the number $\Upsilon(m,n;q)$ is significantly larger than the number of traditional LFSRs of order $n$ over $\mathbb F_{q^m}$, namely, $\phi(q^{mn}-1)/n$, and this is partly a reason why $\sigma$-LFSRs are deemed superior than the LFSRs. 

\begin{remark}
The significance of the power of $q$ in $\Upsilon(m,n;q)$ is not completely clear. We merely mention that $q^{m(m-1)}$ 
is the number of nilpotent $m\times m$ matrices over $\Fq$, thanks to an old result of Fine and Herstein \cite{FH} 
(see \cite{Crabb} or \cite{Ger} for a more accessible proof). Consequently, $|GL_m(\mathbb F_q)|q^{m(m-1)(n-1)}$ is the number 
of $n$-tuples $(C_0, C_1, \dots , C_{n-1})$ of $m\times m$ matrices over $\Fq$ where $C_0$ is nonsingular and $C_1, \dots , C_{n-1}$ are nilpotent. However, the relation of such tuples with primitive $\sigma$-LFSRs is not at all clear. 
\end{remark}

\section{Block Companion Matrices}
\label{bcm}

By a $(m,n)$-\emph{block companion matrix} over $\Fq$ we mean $T\in M_{mn}(\Fq)$ of the form
\begin{equation}
\label{typeT} 
T =
\begin {pmatrix}
\mathbf{0} & \mathbf{0} & \mathbf{0} & . & . & \mathbf{0} & \mathbf{0} & C_0\\
I_m & \mathbf{0} & \mathbf{0} & . & . & \mathbf{0} & \mathbf{0} & C_1\\
. & . & . & . & . & . & . & .\\
. & . & . & . & . & . & . & .\\
\mathbf{0} & \mathbf{0} & \mathbf{0} & . & . & I_m & \mathbf{0} & C_{n-2}\\
\mathbf{0} & \mathbf{0} & \mathbf{0} & . & . & \mathbf{0} & I_m & C_{n-1}
\end {pmatrix}, 
\end{equation}
where $C_0, C_1, \dots , C_{n-1}\in M_m(\Fq)$ and $I_m$ denotes the $m\times m$ identity matrix over $\Fq$, while $\mathbf{0}$ indicates the zero matrix in $M_m(\Fq)$. The set of all $(m,n)$-block companion matrices over $\Fq$ shall be denoted by $\bcm$. Using a Laplace expansion or a suitable sequence of elementary column operations, we see that if $T\in \bcm$ is given by \eqref{typeT}, then $\det T = \pm \det C_0$. Consequently, 
\begin{equation}
\label{nonsingT}
T\in \GL_{mn}(\Fq) \Longleftrightarrow C_0\in \GL_m(\Fq).
\end{equation} 
It may be noted that the block companion matrix \eqref{typeT} is the state transition matrix for the $\sigma$-LFSR \eqref{sigmalfsr}.

The following elementary observation reduces the calculation of a $mn\times mn$ determinant to an $m\times m$ determinant. 
It is implicit in \cite{Zeng4} in the binary case, while a proof in the general case can be gleamed from \cite[Thm. 4 and its proof]{N1}. 

\begin{lemma}
\label{mntom}
Let $T\in \bcm$ be given by \eqref{typeT} and let $F(X)\in M_m\left(\Fq[X]\right)$ be defined by 
$F(X) := I_m X^n - C_{n-1}X^{n-1} - \cdots - C_1 X - C_0$. Then the characteristic polynomial of $T$ is equal to 
$\det F(X)$. 
\end{lemma}

As a corollary, we can obtain a more amenable form of Conjecture \ref{conj1}. 

\begin{theorem}
\label{conj2}
Conjecture \ref{conj1} is equivalent to showing that
\begin{equation}
\label{conjform2}
\left|\left\{T\in \bcm \cap \GL_{mn}(\Fq) : o(T) = q^{mn}-1 \right\}\right| = \Upsilon(m,n;q),
\end{equation}
where $\Upsilon(m,n;q)$ is given by the formula \eqref{NoSigmaLFSR} or the equivalent formula \eqref{NoSigmaLFSR2}.
\end{theorem}

\begin{proof}
If $T\in \bcm\cap \GL_{mn}(\Fq)$ is given by \eqref{typeT} and if $F(X)$ is as in Lemma \ref{mntom}, then 
$\det F(X)$ is precisely the polynomial $\Delta (X)$ in Proposition \ref{primsigmalfsr}. Now, the desired result follows readily from Propositions \ref{maxorder} and \ref{primsigmalfsr} together with Lemma~\ref{mntom}. 
\end{proof}

\section{The Characteristic Map}
\label{charmap}

Let 
$$
\bcms : = \left\{T\in \bcm \cap \GL_{mn}(\Fq) : o(T) = q^{mn}-1 \right\}
$$
be the set of Singer cycles among $(m,n)$-block companion matrices, and 
$$
\pmn := \left\{p(X)\in \Fq[X] : p(X) \mbox{ is primitive of degree } mn\right\}
$$
be the set of all primitive polynomials of degree $mn$ over $\Fq$. In view of Proposition~\ref{maxorder}, the restriction to $\bcms$ of the characteristic map
$$
\Phi :M_{mn}(\Fq) \to \Fq[X] \quad \mbox{ defined by } \quad \Phi(T): = \det\left(XI_{mn} - T \right)
$$
gives a map from $\bcms$ to $\pmn$, which we shall denote by $\Psi$. Clearly, 
$$
\bcms = \coprod_{f(X)\in \, \im (\Psi)} \Psi^{-1}\left(f(X)\right) ,
$$
where, as usual, $\coprod$ denotes disjoint union, $\im (\Psi)$ denotes the image of $\Psi$, and $\Psi^{-1}\left(f(X)\right):= \left\{T\in \bcms : \Psi(T) = f(X)\right\}$ denotes the fiber of $f(X)$ for  any $f(X)\in \pmn$. Thus, to prove \eqref{conjform2}, it suffices to determine 
$\im (\Psi)$ and the cardinality of each of the fibers. The former is answered by the following.

\begin{theorem}
\label{surjprop}
The map 	$\Psi: \bcms \to \pmn$ is surjective. 
\end{theorem} 

\begin{proof}
Let $f(X)\in \pmn$. By Proposition \ref{PrimPolyElt}, there is a primitive element $\gamma$ of $\Fqmn^*$ such that $f(\gamma)=0$. Since $f(X)\in \Fq[X]$, the Frobenius automorphism $x\longmapsto x^q$ of $\Fqmn$ permutes the roots of $f(X)$, and thus $\gamma, \gamma^q, \gamma^{q^2}, \dots, \gamma^{q^{mn-1}}$ are precisely the $mn$ distinct roots of 
$f(X)$. Hence 
$$
f(X) = \prod_{j=0}^{m-1} f_j(X)  \quad \mbox{ where } \quad f_j(X) := \prod_{i=0}^{n-1} \left(X - \gamma^{q^{im + j}}\right)\quad \mbox{ for } j=0,\dots , m-1.
$$
Note that the map given by $x\longmapsto x^{q^m}$ is a generator of the Galois group of $\Fqmn$ over $\Fqm$, and for each $j=0,\dots , m-1$, it permutes the roots of $f_j(X)$ among themselves, and so $f_j(X)\in \Fqm[X]$. Moreover, since $q^j$ and $q^{mn}-1$ are relatively prime, we see that each $f_j(X)$ is the minimal polynomial over $\Fqm$ of a primitive element of $\Fqmn^*$, namely, $\gamma^{q^j}$, and thus $f_j(X)$ is a primitive polynomial in $\Fqm [X]$; in particular,  $f_j(X)$ is irreducible in $\Fqm [X]$ and $f_j(0)\ne 0$ for $j=0,\dots , m-1$. Write
$$
f_0(X) = X^n - \beta_{n-1} X^{n-1} - \cdots - \beta_1 X - \beta_0 \quad \mbox{ where }\beta_0, \beta_1, \dots , \beta_{n-1}\in \Fqm.
$$
Let $B= {\sf C}_{f_0}$ be the companion matrix of $f_0(X)$. By the Cayley-Hamilton Theorem, $f_0(B)=0$ and hence 
$f(B)=0$. Now, choose a Singer cycle $A\in \GL_m(\Fq)$ and let $g(X)\in \Fq[X]$ be the minimal polynomial of $A$. 
By  Proposition \ref{maxorder}, we have $g(X) \in \pmone$. 
Moreover, $p(X)\longmapsto p(A)$ defines a $\Fq$-algebra homomorphism of $\Fq[X]$ into $M_m(\Fq)$ and its image is 
the group algebra $\Fq[A]$ of the cyclic subgroup of $\GL_m(\Fq)$ generated by $A$ while its kernel is the ideal of 
$\Fq[X]$ generated by $g(X)$. Since $g(X)$ is irreducible of degree $m$, the residue class ring $\Fq[X]/\left\langle g(X) \right\rangle$ is $\Fq$-isomorphic to $\Fqm$. Thus we obtain a $\Fq$-algebra isomorphism 
$
\theta : \Fqm \to \Fq[A], 
$
which induces a $\Fq$-algebra homomorphism 
$$
\widehat{\theta}: M_n\left(\Fqm\right) \to M_n \left(M_m(\Fq)\right)\simeq M_{mn}(\Fq) \quad \mbox{given by} \quad
\widehat{\theta} \left(\left(\beta_{ij}\right)\right) = \left(\theta\left(\beta_{ij}\right)\right)
$$ 
of the corresponding rings of matrices. It may be noted that since $o(A)= q^m-1$, we have 
$\Fq[A] =\left\{\mathbf{0}, A, A^2, \dots , A^{q^m-1}\right\}$, where $\mathbf{0}$ denotes the zero matrix in $M_m(\Fq)$. Now let $C_i := \theta(\beta_i)$ for $i=0, \dots , n-1$, and
let $T\in \bcm$ be the matrix given by \eqref{typeT} corresponding to these 
$m\times m$ matrices $C_0, C_1, \dots , C_{m-1}$. Note that since $\beta_0 = f_0(0)\ne 0$ and $\theta$ is an isomorphism,  $C_0$ is nonsingular and hence by \eqref{nonsingT}, $T\in \bcm \cap \GL_{mn}(\Fq)$. Also note that 
$T = \widehat{\theta}(B)$. Now, since $f(B)=0$ and $\widehat{\theta}$ is a $\Fq$-algebra homomorphism, it follows that $f(T)=0$. Moreover, since  $f(X)\in \Fq[X]$ is primitive of degree $mn$, it must be the minimal polynomial of $T$ and further, by Proposition \ref{maxorder}, we see that $o(T)= q^{mn}-1$ and $f(X)$ is the characteristic polynomial of $T$. Thus, $T\in \bcms$ and $\Phi (T) = f(X)$. This proves that $\Psi$ is surjective. 
\end{proof}

As for the fibers of $\Psi$, we propose the following.

\begin{conjecture}[Fiber Conjecture]
\label{twoconj}
For any $f(X)\in \pmn$, the cardinality of the fiber  $\Psi^{-1}\left(f(X)\right):= \left\{T\in \bcms : \Psi(T) = f(X)\right\}$ is independent of the choice of $f(X)$ and, in fact, given by the following formula:
	$$
\left|\Psi^{-1}\left(f(X)\right)\right| =q^{m(m-1)(n-1)} \displaystyle \prod_{i=1}^{m-1}(q^m-q^i). 	
$$
\end{conjecture}

It is clear that Conjecture \ref{twoconj} together with Theorem \ref{surjprop} implies Conjecture \ref{conj1}. 
We remark 
that the fibers of the ambient map $\Phi$ have been studied in the literature (cf. \cite{Ger,Reiner}). 
The Fiber Conjecture facilitates a connection between Conjecture \ref{conj1} 
and a question of Niederreiter (which is still open) as indicated below. 

\begin{remark}
\label{rem:Nied}
Let $\alpha$ be a primitive element of $\F^*_{q^{mn}}$. A subspace $W$ of $\F_{q^{mn}}$ of dimension $m$ is said to be 
{\em $\alpha$-splitting} if $\F_{q^{mn}} = W \oplus \alpha W \oplus \cdots \oplus \alpha^{n-1}W$. Niederreiter [\ref{Niederreiter2}, p. 11] asks 
for the total number of $\alpha$-splitting subspaces of dimension $m$.  In view of Proposition~\ref{PrimPolyElt},  fixing a primitive 
element of $\F^*_{q^{mn}}$ is essentially the same as fixing a primitive polynomial in $\Fq[X]$ of degree $mn$. Now let 
$L_{\alpha} : \F_{q^{mn}} \to \F_{q^{mn}}$ be the linear transformation defined by $L_{\alpha}(x) := \alpha x$. Note that the characteristic polynomial of $L_{\alpha}$ is precisely the minimal polynomial of $\alpha$. 
Moreover, if a subspace $W$ of dimension $m$ is $\alpha$-splitting and $\{u_1, \dots , u_m\}$ is an ordered basis of $W$, then ${\mathcal B}^{\alpha}_{(u_1, \dots , u_m)}= \{ u_1, \dots , u_m, \alpha u_1, \dots , \alpha u_m, \dots , \alpha^{n-1}u_1, \dots , \alpha^{n-1}u_m\}$ is a $\Fq$-basis of $\F_{q^{mn}}$ and with respect to this ordered basis, the matrix of $L_{\alpha}$ is a $(m,n)$-block companion matrix. Moreover, thanks to Proposition \ref{maxorder}, this block companion matrix is  a Singer cycle. 
Conversely, a Singer cycle in $\GL_{mn}(\Fq)$ of the form \eqref{typeT} 
must be the matrix of $L_{\alpha}$ with respect to a basis of the form ${\mathcal B}^{\alpha}_{(u_1, \dots , u_m)}$ and then 
$\{u_1, \dots , u_m\}$ clearly spans a $\alpha$-splitting subspace.  In this way, the enumeration of $\alpha$-splitting subspaces of dimension $m$ is essentially equivalent to the determination of cardinalities of the fibers of 
$\Psi$. We refer to the forthcoming paper \cite{GS} for more on this equivalence and some further progress on Conjectures \ref{conj1} and \ref{twoconj}. 
\end{remark}

\section{The Case $n=1$}
\label{proofs}

As noted in the introduction, when $m=1$, \eqref{NoSigmaLFSR} reduces to \eqref{NoLFSR} and hence Conjecture \ref{conj1} readily follows from Proposition \ref{PrimPolyElt}. Also, when $m=1$, the map $\Psi$ is clearly bijective and hence 
Conjecture \ref{twoconj} holds trivially. We will show below that when $n=1$, both the conjectures follow from the structure of Singer cycles. 

\begin{theorem}
\label{onecase}
If $n=1$, then Conjecture \ref{conj1} as well as  Conjecture \ref{twoconj} hold in the affirmative. 
\end{theorem}

\begin{proof}
Suppose $n=1$. Then $\bcms$ is simply the set of all Singer cycles in $\GL_m(\Fq)$.  By Proposition \ref{conjugateSinger}, $\GL_m(\Fq)$ acts transitively on the set of all Singer subgroups by conjugation, and 
hence the number of Singer subgroups of $\GL_m(\Fq)$ is given by $\left|\GL_m(\Fq)\right|/|N|$, where $N$ denotes
the normalizer of a Singer subgroup of $\GL_m(\Fq)$. Moreover, by Proposition \ref{normalizer}, we see that 
$|N| = m (q^m-1)$. 
Finally, since any Singer subgroup of $\GL_m(\Fq)$ contains $\phi(q^m-1)$ generators, i.e., $\phi(q^m-1)$ Singer cycles, it follows that 
$$
\left|\bcmsone\right| = 
\frac{\left|\GL_m(\Fq)\right|}{m (q^m-1)} \phi(q^m-1) = \Upsilon (m, 1; q).
$$ 
Thus, in view of Theorem \ref{conj2}, Conjecture \ref{conj1} is established when $n=1$. To show more generally,
that Conjecture \ref{twoconj} holds in the affirmative when $n=1$, let 
$f(X)\in \pmone$ and $T\in \bcmsone$ be such that $\Psi(T) = f(X)$. By Proposition \ref{maxorder}, the minimal polynomial as well as the characteristic polynomial of $T$ is $f(X)$. In particular, $T$ and the companion matrix $\Cf$ of $f(X)$ have the same set of invariant factors, and therefore, they are similar (cf. [\ref {bib2'}, p. VII.32]). It follows that  $\Psi^{-1}(f(X)) = \{P^{-1}\Cf P : P \in \GL_m(\Fq)\}$. Consequently, 
$$
\left|\Psi^{-1}(f(X))\right| = \frac{\left|\GL_m(\Fq)\right|}{\left|Z(\Cf)\right|} \quad \mbox{where} \quad
Z(\Cf):= \left\{P \in \GL_m(\Fq) : \Cf P = P\Cf \right\}.
$$
Further, the linear transformation of $\Fqm\simeq \Fq^m$ corresponding to $\Cf$ is cyclic and hence by a theorem of Frobenius \cite[Thm. 3.16 and its Corollary]{jac}, the centralizer $Z(\Cf)$ of $\Cf$ consists only of polynomials in $\Cf$.  
Now, the $\Fq$-algebra of polynomials in $\Cf$ is readily seen to be isomorphic to $\Fq[X]/\left\langle f(X)\right\rangle$, and so its cardinality is $q^m$. Consequently, $Z(\Cf) = \{\Cf^j : 0  \le j < q^m\}$ and 
$|Z(\Cf) |= q^m-1$. Thus, 
$$
\left|\Psi^{-1}(f(X))\right| = \frac{\left|\GL_m(\Fq)\right|}{q^m-1} =  \displaystyle \prod_{i=1}^{m-1}(q^m-q^i),
$$
as desired. 
\end{proof}

\begin{remark}
An alternative proof of Conjecture \ref{twoconj} in the case $n=1$ can be obtained using the Reiner-Gerstenhaber formula for the number of square matrices over $\Fq$ with the given characteristic polynomial (cf. \cite[Thm. 2]{Reiner} and \cite[\S 2]{Ger}) together with Proposition \ref{maxorder}. 
\end{remark}
\section{Examples} 
\label{examples}

In this section we outline some small examples to illustrate Conjecture \ref{conj1} and its refined version Conjecture \ref{twoconj}. Throughout, we take $q=2$ and for $1\le i,j\le 2$, we let $e_{ij}$ denote the $2\times 2$ matrix 
over $\Fq$ with $1$ in $(i,j)^{\rm th}$ place and $0$ elsewhere.  Also, let $I=e_{11}+e_{22}$ be the $2\times 2$ identity matrix and $J= e_{11}+e_{12}+e_{21}+e_{22}$ be the 
$2\times 2$ matrix with all the entries equal to $1$.

\begin{example}
\label{exa1}
Consider $m=2$ and $n=2$. There are only $2$ primitive polynomials of degree $2\times 2 = 4$ over $\mathbb F_2$ and, in fact, we have ${\EuScript P}(4,2) = \{x^4+x+1, x^4+x^3+1\}$. 
It is easily verified that $|{\rm BCMS}(2,2;2)|=16$, i.e., the number of nonsingular $(2,2)$-block companion matrices over $\mathbb F_2$ of order $2^4 -1 = 15$ is $16$, as predicted by Conjecture \ref{conj1}. Moreover, the elements 
$$
T= \begin{pmatrix} \mathbf{0} & C_0 \\ I & C_1 \end{pmatrix}
$$
of ${\rm BCMS}(2,2;2)$ for which $\Psi(T) = x^4+x+1$ are precisely those for which the corresponding pair $(C_0, C_1)$ of $2\times 2$ matrices is given by either of the following.
\begin{eqnarray*}
&& \left(J- e_{21}, \; e_{21}\right), \; \left(J- e_{21}, \; J\right), \;
\left(e_{12}+ e_{21}, \;  e_{21}\right), \; \left(e_{12}+ e_{21}, \;  e_{12}\right), \\
&& \left(J- e_{11},  \; I \right), \; \left(J- e_{22}, \; I\right), \; 
\left(J- e_{12},  \;  e_{12}\right), \; \left(J- e_{12}, \;  J\right).
\end{eqnarray*}
On the other hand, $T\in {\rm BCMS}(2,2;2)$ for which $\Psi(T) = x^4+x^3+1$ are precisely those for which the corresponding pair $(C_0, C_1)$ 
is given by either of the following.
\begin{center}
$\left(J- e_{21}, \; e_{21}+e_{22}\right), \; \left(J- e_{21}, \; e_{11}+ e_{21}\right), \;
\left(e_{12}+ e_{21}, \;  e_{22}\right), \; \left(e_{12}+ e_{21}, \;  e_{11}\right)$, \\
$\left(J- e_{11},  \; J- e_{11} \right), \; \left(J- e_{22}, \; J- e_{22}\right), \; 
\left(J- e_{12},  \; e_{12}+e_{22}\right), \; \left(J- e_{12}, \;  e_{11}+e_{12}\right)$.
\end{center}

\noindent
Thus, both the fibers have cardinality $8$, as predicted by Conjecture \ref{twoconj}.
\end{example}

\begin{example}
\label{exa2}
Consider $m=2$ and $n=3$. Then ${\EuScript P}(6,2)$ consists of six polynomials, namely,  $x^6+x^5+x^4+x+1$, $x^6+x+1$, $x^6+x^5+x^3+x^2+1$, $x^6+x^5+1$, $x^6+x^4+x^3+x+1$, and $x^6+x^5+x^2+x+1$. The fibers of $\Psi$ for each of these consists of $32$ elements of ${\rm BCMS}(2,3;2)$, which together, constitute the 192 elements of ${\rm BCMS}(2,3;2)$. 
It is seen, therefore,  that Conjecture \ref{conj1} as well as Conjecture \ref{twoconj} is valid in this case. 
\end{example}

\section*{Acknowledgments}
  
We are grateful to Surinder Singh Bedi, Gilles Lachaud, Harish Pillai, Samrith
Ram, Sivaramakrishnan Sivasubramanian, and Patrick Sol\'e for helpful discussions, and to Harald Niederreiter for helpful correspondence.
The last two authors are also grateful to
Director, SAG for his permission to publish this paper.

\end{document}